\documentclass[reqno]{amsart}
\usepackage{amssymb,amsmath,amsthm,amsfonts,fancyhdr}
\usepackage{tikz}
\usepackage{indentfirst}
\usepackage{mathrsfs}
\usepackage{setspace}
\usepackage{color}
\usepackage{cite}
\usepackage{needspace}
\allowdisplaybreaks[2]
\usepackage{bm}
\usepackage{geometry}

\numberwithin{equation}{section}

\newtheorem{thm}{Theorem}[section]
\newtheorem{defn}[thm]{{Definition}}
\newtheorem{lemma}[thm]{Lemma}
\newtheorem{cor}[thm]{Corollary}
\newtheorem{re}[thm]{Remark}
\DeclareMathOperator{\osc}{osc}

  \numberwithin{equation}{section}
  \numberwithin{figure}{section}

\begin{document}

\title[Boundary Regularity for  Fully Nonlinear Degenerate Elliptic Equations]{Boundary Regularity for Fully Nonlinear Degenerate  Elliptic  Equations}

\author{Dongsheng Li}
\author[Xuemei Li]{Xuemei Li$^{*}$}
\thanks{$^{*}$Corresponding author: Xuemei Li. Email: \texttt{sms\_lixm@ujn.edu.cn}.}
\author{Lihe Wang}

\address{D. Li: School of Mathematics and Statistics, Xi'an Jiaotong University, Xi'an, P.R.China 710049}
\address{X. Li: School of Mathematical Sciences, University of Jinan, Jinan, P.R.China 250022}
\address{L. Wang: Department of Mathematics, University of Iowa, Iowa City, USA 52242}

\email{lidsh@mail.xjtu.edu.cn}
\email{sms\_lixm@ujn.edu.cn}
\email{lihe-wang@uiowa.edu}

\begin{abstract}
This paper investigates fully nonlinear degenerate elliptic equations of the form $|Du|^\gamma F(D^2u) = f$ with $\gamma>0$ {and $F(0)=0$}.
	Using a unified compactness and perturbation framework, we {establish} $C^1$ regularity on $C^{1,\mathrm{Dini}}$ domains and {$C^{1,\alpha_0}$ regularity on $C^{1,\alpha}$ domains}, {where $\alpha_0=\min\{\alpha,1/(1+\gamma)\}$}. {In addition}, we construct a counterexample to show that {$W^{2,\delta}$ regularity does not hold for any $\delta>0$}.

\end{abstract}



\maketitle
\enlargethispage{4pt}

\section{Introduction}
In this paper, we study the fully nonlinear degenerate elliptic equation
\begin{equation}\label{1.1}
    |Du|^\gamma F(D^2u)=f
    \quad \text{in } \Omega,
\end{equation}
where $\gamma>0$ and $\Omega\subset\mathbb{R}^n$ is a domain.
Throughout the paper, we assume that $F(0)=0$ and that $F$ is
uniformly elliptic: there exist constants $0<\lambda\leq\Lambda$
such that
\begin{equation}\label{ue}
    \lambda\|N\|
    \leq F(M+N)-F(M)
    \leq \Lambda\|N\|,
    \qquad M,N\in\mathbb{S}(n),\quad N\geq0,
\end{equation}
where $\mathbb{S}(n)$ denotes the space of real symmetric
$n\times n$ matrices and $\|\cdot\|$ is the operator norm.
Although $F$ is uniformly elliptic, equation~\eqref{1.1}
degenerates at points where $Du=0$.
Combining compactness and perturbation arguments, we establish
boundary pointwise $C^1$ regularity on $C^{1,\mathrm{Dini}}$ domains
and $C^{1,\alpha_0}$ regularity on $C^{1,\alpha}$ domains, where
$\alpha_0=\min\{\alpha,1/(1+\gamma)\}$.
We also construct a counterexample showing that, for every
$\delta>0$, bounds on $u$ and $f$ in $L^\infty$ alone do not yield
a uniform interior $W^{2,\delta}$ estimate.

When $\gamma=0$, \eqref{1.1} is uniformly elliptic. Caffarelli\cite{caffarelli1989} established the
Alexandrov--Bakelman--Pucci estimate and the Harnack inequality
for viscosity solutions, together with interior $C^\alpha$,
$C^{1,\alpha}$, $C^{2,\alpha}$, and $W^{2,p}$ regularity under
the respective structural and data assumptions. We refer to
\cite{caffarelli1995} for a detailed account.
Winter \cite{winter2009} extended Caffarelli's interior $W^{2,p}$
estimates to the boundary, first for flat boundaries and then
for $C^{1,1}$ domains by a boundary flattening argument.
Sirakov \cite{sirakov2010} obtained boundary $C^\alpha$ estimates
on domains satisfying a uniform exterior cone condition.
Subsequently, Silvestre and Sirakov \cite{silvestre2014}
established boundary $C^{1,\alpha}$ and $C^{2,\alpha}$
regularity on $C^2$ and $C^{2,\tilde{\alpha}}$ domains,
respectively, for suitable H\"older exponents, using
boundary flattening.
Using perturbation and compactness arguments,
Lian and Zhang \cite{lian2020} established pointwise boundary
$C^{1,\alpha}$ and $C^{2,\alpha}$ regularity on
$C^{1,\alpha}$ and $C^{2,\alpha}$ domains, respectively.
Combining the above boundary $C^{1,\alpha}$ regularity {with} interior $W^{2,p}$ regularity, Li, Li and Zhang \cite{li2023} established boundary $W^{2,p_0}$ regularity on $C^{1,\alpha}$ domains for $p_0<\min\{1/(1-\alpha),p/n\}$ {using a} Whitney decomposition.

When $\gamma>0$, \eqref{1.1} degenerates as $|Du|$ tends to zero.
D\'avila, Felmer and Quaas \cite{davila2009} established the
Alexandrov--Bakelman--Pucci estimate, and Imbert
\cite{imbert2011} proved the Harnack inequality.
Imbert and Silvestre \cite{imbert2013} obtained interior
$C^{1,\alpha}$ regularity, which Birindelli and Demengel
\cite{birindelli2014} extended to the boundary for general
viscosity solutions on $C^2$ domains.
Subsequently, Ara\'ujo and Sirakov \cite{araujo2023}
established sharp boundary $C^{1,\alpha}$ regularity on
$C^2$ domains, with an explicit optimal exponent determined
by the degeneracy and the regularity of the boundary data.
In the parabolic setting, Yun \cite{yun2025} obtained
boundary $C^{1,\alpha}$ estimates for analogous degenerate
fully nonlinear equations with convex or concave operators,
under suitable intrinsic regularity assumptions on the
parabolic boundary and the boundary data, using regularization
and compactness methods.
For results on other boundary value problems, we refer to
\cite{banerjee2022,birindelli2022,ricarte2020,byun2025}.

This paper establishes pointwise boundary regularity for
\eqref{1.1} with $\gamma>0$.
Theorems~\ref{t1.1}--\ref{t1.2} give boundary $C^1$ regularity
under pointwise $C^{1,\mathrm{Dini}}$ assumptions on the boundary
and the Dirichlet data and  boundary
$C^{1,\alpha_0}$ regularity with
$\alpha_0=\min\{\alpha,1/(1+\gamma)\}$ under pointwise $C^{1,\alpha}$ assumptions, respectively.
For \eqref{1.1}, these results weaken the $C^2$ boundary
assumption in \cite{birindelli2014,araujo2023}.
In particular, the H\"older result retains the sharp boundary
exponent obtained in \cite{araujo2023}.
We also construct a counterexample to show that {$W^{2,\delta}$ regularity does not hold for any $\delta>0$}.
More precisely, for every
$\delta>0$, bounds on $u$ and $f$ in $L^\infty$ alone do not yield
a uniform interior $W^{2,\delta}$ estimate.

{We now} state our main results. The first two {concern} boundary pointwise regularity, while the last {shows that} {$W^{2,\delta}$ regularity does not hold for any $\delta>0$}{.}

The first result {concerns} \emph{boundary pointwise $C^{1}$ regularity} on $C^{1,\mathrm{Dini}}$ domains ({see Definitions} 2.2 and {2}.3 for the {relevant} notions and {notation}).

\begin{thm}[\textbf{Boundary $C^{1}$ regularity on   $C^{1,\mathrm{Dini}}$ domains}]\label{t1.1}
	Let $\gamma>0$ and {let} $F$ be a fully nonlinear uniformly elliptic operator satisfying \eqref{ue} {and $F(0)=0$}. Suppose that $\Omega$ satisfies the $C^{1,\mathrm{Dini}}$ condition at $0\in \partial \Omega$ with Dini {modulus of continuity} $\omega_\Omega$ and {that} $u \in C(\bar\Omega)$ satisfies
	\begin{equation}\label{equ}
		\left\{
		\begin{aligned}
			|Du|^{\gamma}F(D^2 u) &=f   &\mathrm{in} \ \   &\ \Omega\cap B_1;\\
			u&=g   &\mathrm{on}  \ \  &\partial\Omega\cap B_1,
		\end{aligned}
		\right.
	\end{equation}
	where
	\begin{equation}\label{fg}
		f\in C(\Omega_1)\cap L^{\infty}(\Omega\cap B_1)\ \ \ \mathrm{and} \ \ \ g\in C^{1,\mathrm{Dini}}(0) 
	\end{equation}	
with Dini {modulus of continuity} $\omega_g$. Then $u\in C^{1}(0)${; that is}, {there} exist an affine function $L$ and a modulus of continuity $\sigma$ such that
	\begin{equation}\label{l12}
	|u(x) - L(x)| \leq C |x|\sigma(|x|) \left( \|u\|_{L^\infty(\Omega \cap B_1)} + \|f\|_{L^\infty(\Omega \cap B_1)}^{\frac{1}{1+\gamma}} + \|g\|_{C^{1,\mathrm{Dini}}(0)} \right), \quad \forall x \in \Omega \cap B_{r_0}, 
\end{equation}
	and
	\begin{equation}\label{l13}
	|Du(0)|=|DL| \leq C \left( \|u\|_{L^\infty(\Omega \cap B_1)} + \|f\|^{\frac{1}{1+\gamma}}_{L^\infty(\Omega \cap B_1)} + \|g\|_{C^{1,\mathrm{Dini}}(0)} \right), 
\end{equation}
	where {$C$, $r_0$ and $\sigma$} depend on $n$, $\lambda$, $\Lambda$, $\gamma$, $\omega_g$, $\omega_\Omega$ {and $\|\partial\Omega\|_{C^{1,\mathrm{Dini}}(0)}$}.
\end{thm}

The second result {concerns} \emph{boundary pointwise $C^{1,{\alpha_0}}$ regularity} on $C^{1,\alpha}$ domains ({see Definitions} 2.2 and {2}.3 for the {relevant} notions and {notation}).

\begin{thm}[\textbf{Boundary {$C^{1,\alpha_0}$ regularity on $C^{1,\alpha}$ domains}}]\label{t1.2}
	Let $\gamma>0$ {and $0<\alpha\le1$, and let $F$ satisfy \eqref{ue} and $F(0)=0$}. Suppose that $\Omega$ satisfies {the} $C^{1,\alpha}$ condition at $0\in \partial \Omega$ and {that} $u \in C(\bar\Omega)$ satisfies
	\eqref{equ}
	with
	\begin{equation}\label{fg-holder}
		f\in C(\Omega\cap B_1)\cap L^{\infty}(\Omega\cap B_1)\ \ \ \mathrm{and} \ \ \ g\in C^{1,\alpha}(0).
	\end{equation}	
Set 	\begin{equation}\label{index}
	\alpha_0=\min\{\alpha,1/(1+\gamma)\}.
\end{equation}	
	Then $u\in C^{1,\alpha_0}(0)${; that is}, there exists an affine function $l$ such that
	\begin{equation}\label{l1}
		|u(x)-l(x)|\leq C|x|^{1+\alpha_0}
		\left(||u||_{L^\infty(\Omega\cap B_1)}+||f||_{L^\infty(\Omega\cap B_1)}^{\frac{1}{1+\gamma}}+||g||_{C^{1,\alpha}(0)}\right),\ \forall x\in \Omega\cap B_{r_1}
	\end{equation}
	and
	\begin{equation}\label{l2}
		|Du(0)|=|Dl|\leq C\left(||u||_{L^\infty(\Omega\cap B_1)}+||f||_{L^\infty(\Omega\cap B_1)}^{\frac{1}{1+\gamma}}+||g||_{C^{1,\alpha}(0)}\right),
	\end{equation}
	where $C$ and $r_1$ depend on $n$, $\lambda$,  $\Lambda$, $\gamma$, $\alpha$, $\alpha_0$ and $||\partial\Omega||_{C^{1,\alpha}(0)}$.
\end{thm}

The last result states that
$W^{2,\delta}$ regularity does not hold for any $\delta>0$.
\begin{thm}\label{t1.4}
	Let $\gamma>0$. Suppose that $u \in L^\infty(\Omega)\cap C(\Omega)$ satisfies
$$|Du|^{\gamma}F(D^2 u)  =f   \ \ \mathrm{in}  \ \Omega,$$
where $f\in  L^\infty(\Omega).$
However, $W^{2,\delta}$ regularity does not hold for any $\delta>0$.
\end{thm}

To obtain boundary $C^1$ and $C^{1,\alpha_0}$ estimates
for \eqref{1.1} with $\gamma>0$, we construct affine
approximations with errors bounded by $Cr\sigma(r)$
and $Cr^{1+\alpha_0}$, respectively, in balls of radius $r$
centered at the boundary point.
Here, $\sigma$ is a modulus of continuity.
A key difficulty arises when we subtract affine functions
during the iteration.
The resulting equations take the form
$|Du+p|^\gamma F(D^2u)=f$.
The vector $p$ changes with the scale and may become unbounded,
so the estimates must be uniform in $p$.
Imbert and Silvestre \cite{imbert2013} obtained such
continuity estimates in the interior and used them to establish
interior $C^{1,\alpha}$ regularity.
At the boundary, the geometry and the boundary data introduce
additional errors at each scale.
Their interior argument therefore requires further estimates
to handle the pointwise boundary geometry considered here.
Our key new ingredient is a barrier construction near an almost flat boundary.
It gives boundary estimates independent of $p$ (Lemma~\ref{d1}).
We combine these estimates with the interior continuity estimates of
Imbert and Silvestre \cite{imbert2013} (Lemma~\ref{l3.3}).
This provides the boundary control needed for compactness.
We work directly with the curved boundary and treat it as a perturbation
of a flat boundary, so no boundary flattening is needed.

We outline our approach as follows.
After normalizing the data, we first prove an affine approximation
at one fixed scale (Lemma~\ref{k1}).
A compactness argument gives a homogeneous uniformly elliptic limit equation
with zero data on a flat boundary.
Boundary regularity for this limit gives an affine approximation
with an error of order $r^2$.
This yields the required approximation for the original problem
at a smaller scale.
We then rescale and repeat the argument to construct affine functions $L_k$
such that
\begin{equation}\label{e3}
    |u(x)-L_k(x)|\leq \eta^k\sigma(\eta^k),
    \qquad x\in\Omega_{\eta^k},\quad k\geq1,
\end{equation}
where $0<\eta<1$ and $\sigma$ is a modulus of continuity.
For Theorem~\ref{t1.1}, the Dini condition makes the slope increments
summable and ensures convergence of the affine approximations.
The resulting error is $Cr\sigma(r)$ with $\sigma(r)\to0$,
which gives the boundary $C^1$ estimates \eqref{l12}--\eqref{l13}.
For Theorem~\ref{t1.2}, the same argument gives an error of
$Cr^{1+\alpha_0}$.
The bound $\alpha_0\leq\alpha$ controls the boundary errors,
and $\alpha_0\leq1/(1+\gamma)$ keeps the right-hand side small
at every scale.
This proves the boundary $C^{1,\alpha_0}$ estimates
\eqref{l1}--\eqref{l2}.
Finally, Theorem~\ref{t1.4} follows from a family of smooth oscillating solutions.
The solutions and their right-hand sides are uniformly bounded,
but their Hessians are unbounded in $L^\delta$ for every $\delta>0$.

The paper is organized as follows. In Section 2, {we introduce} basic notions, {including} pointwise characterizations of functions and domains, viscosity solutions and solution sets{. In} Section 3, we prove {uniform} estimates for $|Du+p|^{\gamma}F(D^2 u)=f$ {that are} independent of $p$ and {recall several} known preliminary estimates{,} including interior $C^{\beta}$ estimates for $|Du+p|^{\gamma}F(D^2 u)=f$ {that are} independent of $p$ and boundary $C^{2,\bar\alpha}$ estimates for $F(D^2 u)=0$ on flat domains. In Section 4, {we use} the estimates {established} in Section 3 {to prove Theorems} \ref{t1.1}-\ref{t1.2} {by} perturbation and compactness methods. In Section 5, we prove Theorem \ref{t1.4} by constructing a counterexample {showing} that $W^{2,\delta}$ regularity does not hold for any $\delta>0$.

We end this section by introducing some {notation}.

\noindent\textbf{Notation.} \\
1. $e_i=(0,...,0,1,...,0)=i^{th}$ standard coordinate vector.\\
2. $x'=(x_1,x_2,...,x_{n-1})$ and $x=(x',x_n).$\\
3. $\mathbb{R}^n_+=\{x\in \mathbb{R}^n:x_n>0\}.$\\
4. $B_r(x_0)=\{x\in \mathbb{R}^n: |x-x_0|<r\}$ and $B^+_r (x_0) = B_r(x_0) \cap  \mathbb{R}^n_+$.\\
5. $B_r'=\{x'\in \mathbb{R}^{n-1}: |x'|<r\}$.\\
6. $T_r= B_r\cap\{x_n=0\}.$\\
7. $\Omega_r (x_0)= \Omega\cap B_r(x_0)$ and $(\partial\Omega)_r (x_0)= \partial\Omega\cap B_r(x_0)$. We omit $x_0$ when $x_0=0.$\\
8. $\mathbb{S}(n)$ is the space of real $n\times n$ symmetric matrices. {For $M\in\mathbb S(n)$, $\|M\|$ denotes its operator norm.}

\section{Basic Definitions}
\subsection{Pointwise regularity of functions and domains}
\begin{defn}
	A function $\omega: (0,+\infty) \to [0,+\infty)$ is called a Dini function if $\omega$ is {nonnegative and} nondecreasing {and} satisfies the following Dini condition for some $r_0 > 0${:}
\begin{equation}\label{2.0}
I_\omega:=\int_{0}^{r_0} \frac{\omega(r)}{r} dr < \infty.  
\end{equation}
\end{defn}

{We now define} the {geometric} conditions {under} which we prove our main results. Usually, {saying} that $\partial\Omega$ is $C^{k,\alpha}$ near $x_0\in \partial\Omega$ {means} that $\partial\Omega\cap B_{r_0}(x_0)$ (for some $r_0>0$) can be represented as {the} graph of a $C^{k,\alpha}$ function. Here, we use a more general pointwise {definition}. See {Definition} 2.5 in \cite{lian2020pointwise} for {further} details.

\begin{defn}[\textbf{Pointwise $C^{1,\mathrm{\omega}}$ regularity {of} domains}]
	Let $\Omega$ be a bounded domain, {let} $x_0\in\Gamma\subset \partial\Omega$ and {let} $\omega$ be a modulus of continuity. We say that $\Gamma$ is $C^{1,\mathrm{\omega}}$ at $x_0${,} or {write} $\Gamma\in C^{1,\mathrm{\omega}}(x_0)${,} if there exist $K, r_0 > 0$ and a coordinate system $(x_1,\dots,x_n)$ {isometric to the original one} {in which} $x_0 = 0$ {and}
	\begin{equation}\label{1.2}
{\begin{aligned}
B_{r_0}\cap\{x_n>K|x'|\omega(|x'|)\}&\subset B_{r_0}\cap\Omega,\\
B_{r_0}\cap\{x_n<-K|x'|\omega(|x'|)\}&\subset B_{r_0}\cap\Omega^c.
\end{aligned}}
\end{equation}

In particular,	if $\omega$ is a Dini function, we say that $\Gamma \in C^{1,\mathrm{Dini}}(x_0)$ and define
$$
\|\Gamma\|_{C^{1,\mathrm{Dini}}(x_0)} = \min \left\{ K J_\omega\,\bigg|\, \eqref{1.2}\  {\mathrm{ holds\ with }\ K} \right\},\ 
J_\omega:=I_\omega+\omega(r_0).
$$
Similarly, if $\omega(r)=r^{\alpha}$, we say that $\Gamma \in C^{1,\alpha}(x_0)$ and define
$$
\|\Gamma\|_{C^{1,\alpha}(x_0)} = \min \left\{ K\,\bigg|\, \eqref{1.2}\  {\mathrm{ holds\ with }\ K} \right\}.
$$
\end{defn}

Next, we introduce {definitions describing} the pointwise regularity {of} functions. See {Definition} 2.2 in \cite{lian2020pointwise} for {further} details.

\begin{defn}[\textbf{Pointwise $C^{1,\mathrm{\omega}}$ regularity {of} functions}]
	Let $\Omega \subset \mathbb{R}^n$ be a bounded set ({not necessarily} a domain), {let} $\omega$ be a modulus of continuity and {let} {$f:\Omega\to\mathbb R$}. We say that $f$ is $C^{1,\omega}$ at $x_0 \in \Omega${,} or {write} $f \in C^{1,\omega}(x_0)${,} if there exist a vector $l$ and constants $K, r_0 > 0$ such that
	\begin{equation}\label{1.3}
	|f(x) - f(x_0) - l \cdot (x - x_0)| \leq K |x - x_0|  \omega(|x - x_0|), \quad \forall\, x \in \Omega \cap B_{r_0}(x_0). 
\end{equation}
	{We then} denote $l$ by $D f(x_0)$ and define
	$$ \|f\|_{C^1(x_0)} := |f(x_0)|+|l|, \quad
	[f]_{C^{1,\omega}(x_0)} := \min \left\{ K \omega(r_0) \,\Big|\, \eqref{1.3}\ \mathrm{ holds\ with }\ K \right\}$$
	 and $$
	\|f\|_{C^{1,\omega}(x_0)} := \|f\|_{C^1(x_0)} + [f]_{C^{1,\omega}(x_0)}.
$$
	If $f \in C^{1,\omega}(x)$ for {every} $x \in \Omega$ with the same $\omega, r_0$ and
	\[
	\|f\|_{C^{1,\omega}(\bar{\Omega})} := \sup_{x \in \Omega} \|f\|_{C^{1,\omega}(x)} < +\infty,
	\]
	we say that $f \in C^{1,\omega}(\bar{\Omega})$.	

In particular, if $\omega$ is a Dini function, we say that $f \in C^{1,\mathrm{Dini}}(x_0)$ {and} define
	\[
	[f]_{C^{1,\mathrm{Dini}}(x_0)} = \min \left\{ K J_\omega \,\Big|\, \eqref{1.3}\ \text{ holds\ with\ } K \right\}
	\]
	and
	\[
	\|f\|_{C^{1,\mathrm{Dini}}(x_0)} = \|f\|_{C^1(x_0)} + [f]_{C^{1,\mathrm{Dini}}(x_0)}.
	\]

	Similarly, if $\omega(r) = r^\alpha$ ($0 < \alpha \le 1$), we {define} $f \in C^{1,\alpha}(x_0)$ {by}
	\[
	[f]_{C^{1,\alpha}(x_0)} = \min \left\{ K \,\Big|\, \eqref{1.3}\ \text{ holds\ with\ } K \right\}, \quad
	\|f\|_{C^{1,\alpha}(x_0)} = \|f\|_{C^1(x_0)} + [f]_{C^{1,\alpha}(x_0)}.
	\]
	{We} define $f \in C^{1,\alpha}(\bar{\Omega})$ {analogously}.
	
	If $\omega$ is {merely} a modulus of continuity rather than a Dini function, we may simply {write} $f \in C^1(x_0)$ instead of $f \in C^{1,\omega}(x_0)$.
\end{defn}

\begin{re}
	Since we study {boundary} pointwise regularity, {we assume} throughout this paper, {without loss of generality,} that $0 \in \partial\Omega$ and $r_0 = 1$ in {Definitions} 2.2 and {2}.{3}.
\end{re}

\begin{re} If $\Omega$ is a bounded Lipschitz domain, the definitions of $C^{1,\omega}(\bar{\Omega})$ in {Definition} 2.3 are equivalent to the classical definitions (see \cite{lian2020pointwise}). In this paper, we mainly {consider} the cases {in which} $\Omega$ is a bounded domain or $\Omega$ is part of the boundary of a domain.
\end{re}

 Next, we recall some notions and properties {of} viscosity solutions and {Pucci's} extremal operators $\mathcal{M}^{\pm}$. See \cite{caffarelli1989}
and \cite{caffarelli1995} for more details.

\begin{defn}\label{UE}
	$F$ is {uniformly} elliptic if there exist $0<\lambda\leq\Lambda$ (called ellipticity constants) such that{,} for any $M,N\in \mathbb{S}(n)$ (the space of real $n\times n$ symmetric matrices),
	$$\lambda||N||\leq F(M+N)-F(M)\leq \Lambda||N||,\ \forall N\geq 0.$$
\end{defn}

\begin{defn}\label{V}
	A continuous function $u$ in $\Omega$ is a viscosity subsolution (resp. supersolution) of $F(D^2u)=f(x)$ in $\Omega$ {if} the following condition holds:
	if $x_0\in \Omega$, $\varphi\in C^2(\Omega)$ and $u-\varphi$ has a local maximum (resp. minimum) at $x_0${,} then
	$$F(D^2\varphi(x_0))\geq (resp. \leq) f(x_0).$$
	We say that $u$ is a viscosity solution of $F(D^2u)=f(x)$ {if} it is both a viscosity subsolution and a viscosity supersolution.
\end{defn}

\begin{defn}\label{S}
	Let $0<\lambda\leq\Lambda<\infty$. For $M\in {\mathbb{S}(n)}$ (the space of real $n\times n$ symmetric matrices), {Pucci's} extremal operators $\mathcal{M}^{-}$ and $\mathcal{M}^{+}$ are defined as follows:
	$$\mathcal{M}^{-}(M,\lambda,\Lambda)=\lambda\sum\limits_{e_i>0}e_i+\Lambda\sum\limits_{e_i<0}e_i,$$
	$$\mathcal{M}^{+}(M,\lambda,\Lambda)=\Lambda\sum\limits_{e_i>0}e_i+\lambda\sum\limits_{e_i<0}e_i,$$
	where $e_i=e_i(M)$ are the eigenvalues of {$M$}.
\end{defn}

\begin{defn}
	We denote by $S(\lambda,\Lambda,f)$ the space of continuous functions $u$ in $\Omega$ such that
	$$\mathcal{M}^{-}(D^2u,\lambda,\Lambda)\leq f(x)\leq \mathcal{M}^{+}(D^2u,\lambda,\Lambda)$$
	in the viscosity sense in $\Omega$.
\end{defn}

We {conclude} this section by {defining viscosity subsolutions of} $|Du|^{\gamma}F(D^2 u)  =f$ with $\gamma>0$.

\begin{defn}
Let $\gamma>0$ {and $f\in C(\Omega)$}.
We say that $u\in C(\Omega)$ is a viscosity subsolution (resp. supersolution) of
\[
|Du|^\gamma F(D^2u)=f\quad\text{in }\Omega
\]
if, for every $x_0\in\Omega$ and every $\varphi\in C^2(\Omega)$ such that
$u-\varphi$ has a local maximum (resp. minimum) at $x_0$,
\[
|D\varphi(x_0)|^\gamma F(D^2\varphi(x_0))\ge f(x_0)
\quad(\text{resp. }\le).
\]
{This condition also applies when $D\varphi(x_0)=0$, in which
case it reads $0\ge f(x_0)$ (resp. $0\le f(x_0)$). A viscosity solution is both
a subsolution and a supersolution.}
\end{defn}

\begin{re}
{Since $\gamma>0$, the operator is continuous at zero gradients.
In particular, if $u$ is constant on an open ball, a subsolution satisfies
$0\ge f$ there, and a supersolution satisfies $0\le f$ there.
For $|Du+p|^\gamma F(D^2u)=f$, the same definition is used with
$D\varphi$ replaced by $D\varphi+p$.}
\end{re}

\section{{Preliminary} Estimates}\label{sec3}
\begin{defn}
Given ${0<r\leq1}$, we say {that} the oscillation of $\partial\Omega$ in $B_{r}$ is less than $\delta${,} or write $\mathop{osc}\limits_{B_{r}}\partial\Omega
	\leq\delta$
	if ${B_{r}\cap\{x_n>\delta/2\}\subset \Omega_{r}\subset B_{r}\cap\{x_n>-\delta/2\}}$.
\end{defn}

To carry out {the} compactness arguments for {boundary} regularity, we first prove the following lemma {on} uniform $L^{\infty}$ estimates up to the boundary by constructing barrier functions near the boundary{. These estimates provide} a {form} of equicontinuity up to the boundary.
\begin{lemma}\label{d1}
	Let $\gamma>0$ and $p\in \mathbb{R}^n$. Suppose that $\Omega$ is a {domain}, {$f\in C(\Omega_1)$}, and
 ${u\in C(\overline\Omega\cap B_1)}$ satisfies
	\begin{equation*}
		\left\{
		\begin{aligned}
			|Du+p|^{\gamma}F(D^2 u) &=f   &\mathrm{in} \ \   &\ \Omega_1;\\
			u&=g   &\mathrm{on}  \ \  &(\partial\Omega)_1.
		\end{aligned}
		\right.
	\end{equation*} 
	There exists $\delta_0$ depending on $n,\lambda$ and $\Lambda$ such that{,} if
	$$||u||_{L^\infty(\Omega_1)}\leq 1,\
	{||f||_{L^\infty(\Omega_1)}}\leq\delta,\ ||g||_{L^\infty((\partial\Omega)_1)}\leq\delta\ \mathrm{and}\
	\mathop{osc}\limits_{B_1}\partial\Omega\leq\delta$$
	with $0<\delta<\delta_0$,
	then 
	$$||u||_{L^\infty(\Omega_\delta)}\leq C\delta,$$
	where $C$ depends on $n,\lambda$ and $\Lambda$.
\end{lemma}

\begin{proof}
	Let
	\[
M=2\sqrt{n-1}\left(1+2\sqrt{\frac{(n-1)\Lambda}{\lambda}}\right){,}
\qquad\delta_0=\frac1M.
\]
	{Choose} $\epsilon_0>0$ such that
	\begin{equation}\label{e0}
		4-(2+\epsilon_0)(1+\epsilon_0)(M-1)^{\epsilon_0}>0.
	\end{equation}
	We consider
	$$\begin{aligned}
		\psi=\delta&+4n\left(M+\frac{\delta}{\lambda M}\right)(x_n+\delta)-n\left(M^2+\frac{\delta}{\lambda}\right)x_n(x_n+\delta)\\
		&+\frac{\lambda}{2\Lambda(n-1)}\sum\limits_{i=1}^{n-1}\{(M|x_i|-1)^+\}^{2+\epsilon_0}.
	\end{aligned}$$
	{A straightforward} calculation gives
	$$D_i \psi=\frac{\lambda M(2+\epsilon_0)}{2\Lambda(n-1)}\left\{(M|x_i|-1)^+\right\}^{1+\epsilon_0}\frac{x_i}{|x_i|}, \
	i=1,2,...,n-1$$
{Here $x_i/|x_i|$ is defined to be $0$ when $x_i=0$.}
	{Moreover,}
	$$D_n \psi=4n\left(M+\frac{\delta}{\lambda M}\right)-n\left(M^2+\frac{\delta}{\lambda}\right)(2x_n+\delta).$$
	Set
	$${D=\Omega_1\cap\{x_n<\delta_0\}}.$$
	In $D$, we deduce that
	$$|D_i \psi|\leq\frac{\lambda M(2+\epsilon_0)}{2\Lambda(n-1)}(M-1)^{1+\epsilon_0}, \ i=1,2,...,n-1$$
	and $$
	{0<D_n\psi\leq4n\left(M+\frac{\delta}{\lambda M}\right)
\leq4n\left(M+\frac{1}{\lambda M^2}\right)},$$
	Set 
	$$N=\frac{\lambda M(2+\epsilon_0)}{\Lambda}(M-1)^{1+\epsilon_0}+{8n\left(M+\frac1{\lambda M^2}\right)}.$$
	Then
	$$|D\psi|\leq \frac{\lambda M(2+\epsilon_0)}{2\Lambda}(M-1)^{1+\epsilon_0}+{4n\left(M+\frac1{\lambda M^2}\right)}=N/2\ \ \mathrm{in}\ \ D.$$
	
	{We} prove
	$$u\leq C\delta\ \ \mathrm{in}\ \ \Omega_\delta$$ {by considering} two cases: (i) $|p|>N$ and (ii) $|p|\leq N$.
	
	\textbf{Case (i).} {When} $|p|>N$, we
	show that $\psi$ provides an upper barrier in $D$ with the properties
	\begin{equation}\label{b1}
		|D\psi+p|^{\gamma}F(D^2 \psi)\leq-2\delta<f \ \ \mathrm{in}\ \ D
	\end{equation}
	and
	\begin{equation}\label{b2}
		{\limsup_{D\ni y\to x}(u(y)-\psi(y))\leq0\quad(x\in\partial D)}.
	\end{equation}
	{Once} \eqref{b1} and \eqref{b2} {are established}, {comparison with the smooth strict supersolution $\psi$ gives} $$\psi\geq u\ \ \mathrm{in} \ \  D.$$ {Since} $(M|x_i|-1)^+=0$ in $\Omega_\delta$ and $\Omega_\delta\subset D$, {we obtain}
	\begin{equation}
		u\leq \psi=\delta+4n\left(M+\frac{\delta}{\lambda M}\right)(x_n+\delta)-n\left(M^2+\frac{\delta}{\lambda}\right)x_n(x_n+\delta)\leq C\delta \ \ \mathrm{in}\ \ \Omega_\delta,
	\end{equation}
	where $C$ depends on $n$, $\lambda$ and $\Lambda$.
	
	We now prove \eqref{b1} and \eqref{b2}.
A direct calculation gives
	$$\begin{aligned}
		D_{ii} \psi=\frac{\lambda M^2(2+\epsilon_0)(1+\epsilon_0)}{2\Lambda(n-1)}\left\{(M|x_i|-1)^+\right\}^{\epsilon_0}\leq\frac{\lambda M^2(2+\epsilon_0)(1+\epsilon_0)}{2\Lambda(n-1)}(M-1)^{\epsilon_0},i=1,\ldots,n-1
	\end{aligned}$$
	$$D_{ij}\psi=0, i\neq j, i,j=1,2,...,{n},\ D_{nn} \psi=-2n\left(M^2+\frac{\delta}{\lambda}\right).$$
	Consequently,
	\begin{equation*}
		\begin{aligned}
			\mathcal{M}^+(D^2\psi,\lambda/n,\Lambda)&=\lambda D_{nn}\psi/n+\Lambda\sum\limits_{i=1}^{n-1}D_{ii}\psi\\
			&\leq -2(\lambda M^2+\delta)+\frac{\lambda M^2(2+\epsilon_0)(1+\epsilon_0)}{2}(M-1)^{\epsilon_0}\\
			&\leq -2\delta,
		\end{aligned}
	\end{equation*}
	where {the last inequality follows from} \eqref{e0}.
	{Since} $|p|>N$, $|D\psi|\leq N/2$ and $\gamma\geq0$,  $$|D\psi+p|^{\gamma}\geq(|p|-|D\psi|)^{\gamma}\geq(N/2)^\gamma\geq 1.$$
	{By} uniform ellipticity, we deduce that
	$$|D\psi+p|^{\gamma}F(D^2 \psi)\leq|D\psi+p|^{\gamma}\mathcal{M}^+(D^2 \psi,\lambda/n,\Lambda)\leq -2\delta<f$$
	{which establishes} \eqref{b1}.
	It remains to prove \eqref{b2}. Recall that $\delta_0=\frac{1}{M}$.
	$$\psi\geq n \left(\frac{1}{M}+\delta\right)\left(4\left(M+\frac{\delta}{\lambda M}\right)-\left(M+\frac{\delta}{\lambda M}\right)\right)\geq 1\ \ \mathrm{on} \ \ \partial D\cap\{x_n=\delta_0\}.$$
	{Since} $\mathop{osc}\limits_{B_1}\partial\Omega\leq\delta\leq\delta_0$, we have $-\delta/2\leq x_n\leq \delta_0=1/M$ and {hence,} on $(\partial\Omega)_1$,
	$$\begin{aligned}&4n\left(M+\frac{\delta}{\lambda M}\right)(x_n+\delta)-n\left(M^2+\frac{\delta}{\lambda}\right)x_n(x_n+\delta)\\&
		\geq n(x_n+\delta)\left(4\left(M+\frac{\delta}{\lambda M}\right)-\left(M+\frac{\delta}{\lambda M}\right)\right)\geq 0,
	\end{aligned}$$
	which implies that
	$$\psi\geq \delta\geq g\ \ \mathrm{on} \ \ \partial D\cap{(\partial\Omega)_1}.$$
	When $|x|=1$, {the bound} $|x_n|\leq \delta_0$ {gives}
	$$\max_{1\leq i\leq n-1}|x_i|\geq \frac{1}{2\sqrt{n-1}}$$
	and then
	$$\psi\geq \frac{\lambda}{2\Lambda(n-1)}\left(\frac{M}{2\sqrt{n-1}}-1\right)^2
	=\frac{\lambda}{2\Lambda(n-1)}\left(2\sqrt{\frac{(n-1)\Lambda}{\lambda}}
	\right)^2\geq 1\ \ \mathrm{on} \ \ \partial D\cap\partial B_1.$$
	Therefore, we obtain \eqref{b2}.

	\textbf{Case (ii).} {When} $|p|\leq N$, we consider $\psi_p=\psi+|p|(x_n+\delta)$ and show that $\psi_p$ provides an upper barrier in $D$. Since $\mathop{osc}\limits_{B_1}\partial\Omega\leq\delta$, we {have} $x_n+\delta\geq 0$ in $\Omega_1$ and {hence} $\psi_p\geq \psi${,} which implies {that \eqref{b2} holds with $\psi_p$ in place of $\psi$}.
	Since $x_n\leq \delta_0=1/M$ in $D$, we have
	$$D_n \psi\geq 4n\left(M+\frac{\delta}{\lambda M}\right)-3n\delta_0\left(M^2+\frac{\delta}{\lambda}\right)\geq nM$$ and then $|D \psi_p|\geq D_n \psi_p=D_n \psi+|p|\geq M+|p|.$ It follows that
	$$|D\psi_p+p|^{\gamma}\geq(|D\psi_p|-|p|)^{\gamma}\geq M^{\gamma}\geq 1.$$
	Thus,
	$$|D\psi_p+p|^{\gamma}F(D^2 \psi_p)\leq |D\psi_p+p|^{\gamma}\mathcal{M}^+(D^2\psi_p,{\lambda/n,\Lambda})\leq -2\delta<f\ \ \mathrm{in}\ \ D.$$
	{Consequently}, $u\leq \psi_p$ in $D${,} and {hence}
	$$u\leq \psi_p\leq \psi+ 2N\delta \leq C\delta\ \ \mathrm{in}\ \ \Omega_\delta,$$ where $C$ depends on $n$, $\lambda$ and $\Lambda$.
	
	The {lower bound} $u\geq-C\delta$ in $\Omega_\delta$ {follows similarly,} and we omit {the proof}.
\end{proof}

\begin{cor}\label{d1'}
	For any $0<r<1$ and $\epsilon>0$, there exists $\delta>0$ (depending on $n,\lambda,\Lambda,r$ and $\epsilon$) such that if $\Omega$ is a {domain} and ${u\in C(\overline\Omega\cap B_1)}$ satisfies
	\begin{equation*}
		\left\{
		\begin{aligned}
			|Du+p|^{\gamma}F(D^2 u) &=f   &\mathrm{in} \ \   &\ \Omega_1;\\
			u&=g   &\mathrm{on}  \ \  &(\partial\Omega)_1,
		\end{aligned}
		\right.
	\end{equation*}
	with
	$||u||_{L^\infty(\Omega_1)}\leq 1$,
	$||f||_{L^\infty(\Omega_1)}\leq\delta$, $||g||_{L^\infty((\partial\Omega)_1)}\leq\delta$ and
	$\mathop{osc}\limits_{B_1}\partial\Omega\leq\delta$, then
	$$||u||_{L^\infty(\Omega_{\delta}(x_0))}\leq \epsilon,\ \forall x_0 \in (\partial\Omega)_r.$$
\end{cor}

Next, we {recall} two lemmas {that} will be used in the {subsequent} compactness arguments. {See Lemmas} 3 and {6} in \cite{imbert2013} for their {proofs}.

\begin{lemma}\label{l3.3}
	For {every} ${0<r<1}$, there exist $\beta\in (0,1)$ and $C>0$ {depending} only {on} $n,\lambda,\Lambda,\gamma$
	and $r$ such that {every} viscosity solution of $|p+D u|^\gamma F(D^2 u)=f$ in $B_1$ with $\osc_{B_1} u\leq 1$ and $||f||_{L^\infty(B_1)}\leq\epsilon_0<1$ ({where} $\epsilon_0$ depends on $n,\lambda,\Lambda$ and $\gamma$) satisfies
	$$[u]_{C^\beta(B_r)}\leq C.$$
	In particular, the modulus of continuity of $u$ is controlled independently of $p$.
\end{lemma}

\begin{lemma}\label{l3.4}
	Assume that $u$ is a viscosity solution of $|p+D u|^\gamma F(D^2 u)=0$ in $B_1$.
	Then $u$ is a viscosity solution of $F(D^2 u)=0$ in $B_1$.
\end{lemma}

Finally, since our method treats {a nonsmooth} boundary as a perturbation of a flat boundary, we conclude this section with the following lemma {on} pointwise boundary $C^{2,\bar\alpha}$ estimates for $F(D^2 u) = 0$ on flat boundaries with zero boundary values. See Lemma 4.1 in \cite{silvestre2014} and Lemma 2.2 in \cite{lian2020}.

\begin{lemma}\label{l3.1}
Let $F$ be uniformly elliptic with $F(0)=0${,} and {let} $u$ be a viscosity solution {of} the equation
\begin{equation}\label{u00}
	\begin{cases}
		F(D^2 u) = 0 & \text{in } \ B_1^+, \\
		u = 0 & \text{on }\  T_1.
	\end{cases}
\end{equation}
	Then $u \in C^{2,\bar\alpha}(0)$ and, for some constants $a$ and $b_{in}$ ($1 \le i \le n$), {the polynomial
\[
P(x)=a x_n+\sum_{i=1}^n b_{in}x_i x_n
\]
satisfies}
\[
{|u(x)-P(x)|} \le \bar C |x|^{2+\bar\alpha} \|u\|_{L^\infty(B_1^+)},\quad \forall\, x \in B_{1/2}^+,  
\]
\[
{F(D^2P) = 0}  
\]
and
\[
|a| + {\sum_{i=1}^n|b_{in}|} \le \bar C\|u\|_{L^\infty(B_1^+)},
\]
{where} $\bar\alpha$ and $\bar C$ depend on $n,\lambda$ and $\Lambda$.
\end{lemma}

\begin{re}\label{re2}
{It follows immediately} that if $u$ satisfies \eqref{u00}, then
$$ 
\big| u(x) - a x_n \big| \le \bar C |x|^{2} \|u\|_{L^\infty(B_1^+)},\quad \forall\, x \in B_{1/2}^+.$$
\end{re}

\section{Pointwise Boundary Estimates}\label{sec4}

The {proofs} of Theorem \ref{t1.1} and Theorem \ref{t1.2} {proceed} as follows. Uniform $L^\infty$ boundary estimates {were established} in the previous section. Building on these estimates, we {use} compactness and perturbation techniques to show that the solution can be approximated by {an} affine function {when} the prescribed data are {sufficiently} small ({see} Lemma \ref{k1}){. Our} main theorems {then follow from} a standard scaling argument.

To prove our main theorems, we first {establish} the following key compactness lemma.

\begingroup
\begin{lemma}\label{k1}
Let $\bar C\geq1$ be as in Lemma~\ref{l3.1}
and Remark~\ref{re2}.
For every $0<\alpha<1$ and $\gamma>0$, there exist
$\eta\in(0,1/8)$ and $\delta>0$, depending only on
$n,\lambda,\Lambda,\alpha$, and $\gamma$, such that
$\bar C\eta^{1-\alpha}\leq1/2$ and the following holds.

Let $F$ satisfy \eqref{ue} and $F(0)=0$.
Let $p\in\mathbb R^n$ and let $\Omega$ be a domain.
Suppose that $u\in C(\overline\Omega\cap B_1)$
is a viscosity solution of
\[
\begin{cases}
|Du+p|^\gamma F(D^2u)=f
&\text{in }\Omega_1,\\
u=g
&\text{on }(\partial\Omega)_1,
\end{cases}
\]
where $f$ is continuous.
Assume that
\[
\|u\|_{L^\infty(\Omega_1)}\leq1,
\qquad
\|f\|_{L^\infty(\Omega_1)}\leq\delta,
\qquad
\|g\|_{L^\infty((\partial\Omega)_1)}\leq\delta,
\qquad
\mathop{\mathrm{osc}}_{B_1}\partial\Omega\leq\delta.
\]
Then there exists $a\in\mathbb R$ such that
\[
|a|\leq\bar C,
\qquad
\|u-ax_n\|_{L^\infty(\Omega_\eta)}
\leq\eta^{1+\alpha}.
\]
The scale $\eta$ is chosen first, and $\delta$ is fixed afterwards.
Both choices are uniform over all such operators $F$
and all vectors $p\in\mathbb R^n$.
\end{lemma}

\begin{proof}
Fix $\bar C$ as above and choose $\eta\in(0,1/8)$ such that
\[
\bar C\eta^{1-\alpha}\leq\frac12.
\]
Suppose, for contradiction, that no $\delta>0$ works uniformly
under the assumptions of the lemma.
Taking $\delta=1/k$, we obtain domains $\Omega_k$,
operators $F_k$, vectors $p_k\in\mathbb R^n$,
and functions $u_k,f_k,g_k$.
Each $F_k$ satisfies \eqref{ue} with the same ellipticity
constants $\lambda,\Lambda$ and satisfies $F_k(0)=0$.
Each $u_k\in C(\overline{\Omega_k}\cap B_1)$
is a viscosity solution of
\[
\begin{cases}
|Du_k+p_k|^\gamma F_k(D^2u_k)=f_k
&\text{in }(\Omega_k)_1,\\
u_k=g_k
&\text{on }(\partial\Omega_k)_1,
\end{cases}
\]
where $f_k$ is continuous and
\[
\|u_k\|_{L^\infty((\Omega_k)_1)}\leq1,
\qquad
\|f_k\|_{L^\infty((\Omega_k)_1)}\leq\frac1k,
\qquad
\|g_k\|_{L^\infty((\partial\Omega_k)_1)}\leq\frac1k,
\qquad
\mathop{\mathrm{osc}}_{B_1}\partial\Omega_k\leq\frac1k.
\]
However, for every $k\geq1$,
\begin{equation}\label{3.7}
\|u_k-ax_n\|_{L^\infty((\Omega_k)_\eta)}
>\eta^{1+\alpha}
\qquad\text{for every }a\in\mathbb R
\text{ with }|a|\leq\bar C.
\end{equation}

\par\smallskip\noindent
\textbf{Step 1. We establish compactness of the operators and solutions.}
Uniform ellipticity and $F_k(0)=0$ give
\[
|F_k(M)-F_k(N)|\leq n\Lambda\|M-N\|,
\qquad
|F_k(M)|\leq n\Lambda\|M\|
\quad(M,N\in\mathbb S(n)).
\]
Thus $\{F_k\}$ is uniformly bounded and equicontinuous
on every compact subset of $\mathbb S(n)$.
By the Arzel\`a--Ascoli theorem and a diagonal argument,
we can pass to a subsequence such that
\[
F_k\longrightarrow F_\infty
\quad\text{locally uniformly on }\mathbb S(n),
\qquad
F_\infty(0)=0.
\]
Passing to the limit in \eqref{ue}, we obtain
\[
\lambda\|N\|
\leq F_\infty(M+N)-F_\infty(M)
\leq\Lambda\|N\|
\quad(M,N\in\mathbb S(n),\ N\geq0).
\]

For every compact set $K\subset B_{3/4}^+$,
the domain inclusions ensure that a fixed neighborhood of $K$
lies in $\Omega_k\cap B_1$ for all sufficiently large $k$.
By Lemma~\ref{l3.3}, after passing to a subsequence,
\[
u_k\longrightarrow u_0\quad\text{locally uniformly in }B_{3/4}^+,
\qquad \|u_0\|_{L^\infty(B_{3/4}^+)}\le1.
\]

\par\smallskip\noindent
\textbf{Step 2. We identify the limit equation and its boundary values.}
We show that $u_0$ extends continuously to $T_{1/2}$ and satisfies
\begin{equation}\label{cvb}
\begin{cases}
F_\infty(D^2u_0)=0
&\text{in }B_{1/2}^+,\\
u_0=0
&\text{on }T_{1/2}.
\end{cases}
\end{equation}

 We first prove the interior equation. We follow the stability argument in the proof of
\cite[Lemma~2]{imbert2013}.
Every fixed ball compactly contained in $B_{3/4}^+$
lies in $\Omega_k$ for all sufficiently large $k$.
Since $F_k\to F_\infty$ locally uniformly,
the same argument applies to the present sequence of operators.
After passing to a subsequence, either $p_k\to p_0$
or $|p_k|\to\infty$.
If $p_k\to p_0$, then
\[
|p+p_k|^\gamma F_k(M)
\longrightarrow |p+p_0|^\gamma F_\infty(M),
\qquad
f_k\longrightarrow0.
\]
The operator convergence is locally uniform in $(p,M)$.
Viscosity stability therefore gives
\[
|Du_0+p_0|^\gamma F_\infty(D^2u_0)=0
\quad\text{in}\quad B_{1/2}^+.
\]
By Lemma~\ref{l3.4}, $F_\infty(D^2u_0)=0$ in ${B_{1/2}^+}$.
{
If $|p_k|\to\infty$, we divide the equations by $|p_k|^\gamma$.
The normalized operators converge locally uniformly to $F_\infty(M)$,
and the right-hand sides converge to zero.
Viscosity stability therefore gives
\[
F_\infty(D^2u_0)=0\quad\text{in }B_{1/2}^+.
\]
}

We next prove the zero boundary values.
Fix $\epsilon>0$.
Corollary~\ref{d1'}, with $r=3/4$, gives a radius
$0<\rho<1/8$, independent of $k$, such that
for all sufficiently large $k$,
\[
\sup_{\Omega_k\cap B_\rho(y)}|u_k|\leq\epsilon
\qquad
\text{for every }y\in\partial\Omega_k\cap B_{3/4}.
\]
Let $x\in\Omega_k\cap B_{1/2}$ with $x_n\leq\rho/2$.
For sufficiently large $k$, the slab inclusions give
\[
\operatorname{dist}(x,\partial\Omega_k\cap B_{3/4})
\leq\frac{\rho}{2}+\frac{1}{2k}<\rho.
\]
The preceding boundary estimate therefore implies
\begin{equation}\label{3.9}
\sup_{\Omega_k\cap B_{1/2}\cap\{x_n\leq\rho/2\}}
|u_k|\leq\epsilon
\qquad(k\text{ sufficiently large}).
\end{equation}
Every fixed $x\in B_{1/2}^+$ belongs to $\Omega_k$
for all sufficiently large $k$.
Passing to the limit in \eqref{3.9}, we obtain
\[|u_0|\leq\epsilon\qquad \mathrm{in}\quad
B_{1/2}^+\cap\{x_n\leq\rho/2\}.\]
Since $\epsilon>0$ is arbitrary, $u_0$ extends continuously
to $B_{1/2}^+\cup T_{1/2}$ with $u_0=0$ on $T_{1/2}$.

\par\smallskip\noindent
\textbf{Step 3. We obtain an affine approximation and reach a contradiction.}
Applying Lemma~\ref{l3.1} and Remark~\ref{re2} to
\eqref{cvb}, there exists a constant 
$\bar a$ such that
\[
|\bar a|\leq\bar C,
\qquad
|u_0(x)-\bar a x_n|\leq\bar C|x|^2
\quad\text{for }x\in B_{1/4}^+.
\]
By \eqref{3.7}, for each $k$ we can choose
$x_k\in\Omega_k\cap B_\eta$ such that
\[
|u_k(x_k)-\bar a(x_k)_n|>\eta^{1+\alpha}.
\]
After passing to a subsequence, we may assume that
$x_k\to x_*\in\overline B_\eta$.
The slab inclusions imply $(x_*)_n\geq0$.

If $(x_*)_n>0$, local uniform convergence and the continuity
of $u_0$ give
\[
\begin{aligned}
\eta^{1+\alpha}
\leq \lim_{k\to\infty}
      |u_k(x_k)-\bar a(x_k)_n| 
= |u_0(x_*)-\bar a(x_*)_n| 
\leq \bar C\eta^2
\leq \frac12\eta^{1+\alpha},
\end{aligned}
\]
a contradiction.

If $(x_*)_n=0$, fix $\epsilon>0$ and let $\rho$ be as
in \eqref{3.9}.
Since $(x_k)_n\to0$, that estimate gives
$|u_k(x_k)|\leq\epsilon$ for all sufficiently large $k$.
Hence
\[
\eta^{1+\alpha}
< |u_k(x_k)-\bar a(x_k)_n|
\leq \epsilon+\bar C|(x_k)_n|.
\]
Letting $k\to\infty$ gives $\eta^{1+\alpha}\leq\epsilon$,
which is impossible since $\epsilon>0$ is arbitrary.
This proves the lemma.
\end{proof}
\endgroup

\begin{re}
Lemma~\ref{k1} uses the stated slab inclusions and smallness
assumptions, without requiring $C^{1,\mathrm{Dini}}$
regularity of $\partial\Omega$ or $g$.
The pointwise $C^{1,\mathrm{Dini}}$ and $C^{1,\alpha}$
assumptions are used in the iterations for
Theorems~\ref{t1.1} and~\ref{t1.2}, respectively;
see Remarks~\ref{dini} and~\ref{alpha}.
\end{re}

{Motivated} by Theorem 10.8 of \cite{lian2020pointwise} and Theorem 1.6 of \cite{lian2020}, we prove {Theorems} \ref{t1.1}-\ref{t1.2} {using} induction and scaling arguments.

\bigskip

\noindent\textbf{Proof of Theorem~\ref{t1.1}.}

\par\smallskip\noindent
\textbf{Step 1. We normalize the problem.}
Let $\alpha=1/[2(1+\gamma)]$.
Choose $\eta$ as in the proof of Lemma~\ref{k1},
taking $\eta$ small enough that
\[
\eta^\alpha\leq\frac14,
\qquad
\eta^{1/(1+\gamma)}
\left(1+\frac{1+\gamma}{1-\eta}\right)\leq\frac14.
\]
Then fix the corresponding $0<\delta\leq1$.
Let $K_\Omega,K_g\geq0$ be the optimal constants
in the pointwise domain and data estimates
of Definitions~2.2--2.3, respectively.
Then we have
\[
|x_n|\leq K_\Omega|x|\omega_\Omega(|x|),
\qquad
|g(x)-g(0)-Dg(0)\cdot x|
\leq K_g|x|\omega_g(|x|)
\quad\text{on }(\partial\Omega)_1.
\]
Set
\[
J_\Omega(t)=K_\Omega\left(
\omega_\Omega(t)+\int_0^t\frac{\omega_\Omega(s)}{s}\,ds
\right),
\qquad
J_g(t)=K_g\left(
\omega_g(t)+\int_0^t\frac{\omega_g(s)}{s}\,ds
\right),
\]
and
\[
M=\|u\|_{L^\infty(\Omega_1)}
  +\|f\|_{L^\infty(\Omega_1)}^{1/(1+\gamma)}
  +\|g\|_{C^{1,\mathrm{Dini}}(0)}.
\]
By Definitions~2.2--2.3,
\[
J_\Omega(1)=\|\partial\Omega\|_{C^{1,\mathrm{Dini}}(0)},
\qquad
J_g(1)=[g]_{C^{1,\mathrm{Dini}}(0)}\leq M.
\]
The case $M=0$ is immediate.
Assume $M>0$ and choose $0<r\leq\delta/2$ such that
\[
J_\Omega(r)\leq\frac{(1-\eta)\delta}{32\bar C},
\qquad
J_g(r)\leq\frac{1-\eta}{8}J_g(1).
\]

Define
\[
\widetilde\Omega=\Omega/r,
\qquad
\widetilde u(x)
=\frac{u(rx)-g(0)-rDg(0)\cdot x}{M},
\qquad
p=\frac{rDg(0)}{M}.
\]
Then $\widetilde u$ satisfies
\[
\begin{cases}
|D\widetilde u+p|^\gamma
\widetilde F(D^2\widetilde u)=\widetilde f
&\text{in }\widetilde\Omega_1,\\
\widetilde u=\widetilde g
&\text{on }(\partial\widetilde\Omega)_1,
\end{cases}
\]
where
\[
\widetilde F(X)=\frac{r^2}{M}
F\!\left(\frac{M}{r^2}X\right),
\qquad
\widetilde f(x)
=\frac{r^{\gamma+2}}{M^{\gamma+1}}f(rx),
\qquad
\widetilde g(x)
=\frac{g(rx)-g(0)-rDg(0)\cdot x}{M}.
\]
The operator $\widetilde F$ has the same ellipticity constants
as $F$ and satisfies $\widetilde F(0)=0$.

Define
\begin{equation}\label{o1}
\widetilde\omega_\Omega(t)
=\frac{4\bar C K_\Omega}{\delta}\omega_\Omega(rt),
\qquad
\widetilde\omega_g(t)
=\frac{K_g\omega_g(rt)}{J_g(1)},\qquad
\omega(t)=
\max\left\{
\widetilde\omega_\Omega(t),
\widetilde\omega_g(t),
t^{1/(1+\gamma)}
\right\}.
\end{equation}
{When $J_g(1)=0$, we set $\widetilde\omega_g\equiv0$
and interpret $J_g(r)/J_g(1)$ as zero.}
By the choices of $M$, $r$, and $\omega$,
the normalized quantities satisfy
\begin{equation}\label{a3}
\begin{gathered}
\|\widetilde u\|_{L^\infty(\widetilde\Omega_1)}\leq1,
\quad
\|\widetilde f\|_{L^\infty(\widetilde\Omega_1)}\leq\delta;\qquad
|\widetilde g(x)|
\leq\frac{\delta}{2}|x|\omega(|x|),
\quad
|x_n|
\leq\frac{\delta}{4\bar C}|x|\omega(|x|)
\quad\text{on }(\partial\widetilde\Omega)_1,\\
\omega(1)=1,
\qquad
\omega(\eta)
+\frac{1}{1-\eta}\int_0^\eta\frac{\omega(t)}{t}\,dt
\leq\frac12.
\end{gathered}
\end{equation}
For the boundary data estimate, we used
$rJ_g(1)/M\leq\delta/2$.
The last inequality follows from
\[
\begin{aligned}
\omega(\eta)
+\frac{1}{1-\eta}\int_0^\eta\frac{\omega(t)}{t}\,dt
\leq
\frac{1}{1-\eta}
\left(
\frac{4\bar C}{\delta}J_\Omega(r)
+\frac{J_g(r)}{J_g(1)}
\right)
{+}\eta^{1/(1+\gamma)}
\left(1+\frac{1+\gamma}{1-\eta}\right)
\leq\frac12.
\end{aligned}
\]
We will prove that there exist $a\in\mathbb R$ and a modulus
of continuity $\widetilde\sigma$ such that
\begin{equation}\label{star}
|\widetilde u(x)-ax_n|
\leq C|x|\widetilde\sigma(|x|)
\quad\text{for }x\in\widetilde\Omega_\eta,
\qquad
|a|\leq C.
\end{equation}

\par\smallskip\noindent
\textbf{Step 2. We construct affine approximations by induction.}
Define
\begin{equation}\label{lip3}
A_0=1,
\qquad
A_k=\max\{\omega(\eta^k),\eta^\alpha A_{k-1}\}
\quad(k\geq1).
\end{equation}
We construct a sequence of constants $\{a_k\}_{k\geq0}$ such that
\begin{equation}\label{lip1}
\|\widetilde u-a_kx_n\|_{L^\infty(\widetilde\Omega_{\eta^k})}
\leq\eta^kA_k
\quad(k\geq0),
\end{equation}
and
\begin{equation}\label{lip2}
|a_k-a_{k-1}|
\leq\bar C A_{k-1}
\quad(k\geq1),
\end{equation}
where $\bar C$ is as in Lemma~\ref{l3.1}.

For $k=0$, set $a_0=0$.
Then \eqref{lip1} follows from
$\|\widetilde u\|_{L^\infty(\widetilde\Omega_1)}\leq1$.
Assume that \eqref{lip1}--\eqref{lip2} hold for all indices up to $k$.
We prove them for $k+1$.
Set
\begin{equation}\label{s1}
\Omega_k=\widetilde\Omega/\eta^k,
\qquad
u_k(x)=
\frac{\widetilde u(\eta^kx)-a_k\eta^kx_n}{\eta^kA_k},
\qquad x\in(\Omega_k)_1.
\end{equation}
Then $u_k$ satisfies
\[
\begin{cases}
|Du_k+A_k^{-1}(p+a_ke_n)|^\gamma
\eta^kA_k^{-1}\widetilde F(\eta^{-k}A_kD^2u_k)=f_k
&\text{in }(\Omega_k)_1,\\
u_k=g_k
&\text{on }(\partial\Omega_k)_1,
\end{cases}
\]
where
\[
f_k(x)=\eta^kA_k^{-(1+\gamma)}\widetilde f(\eta^kx),
\qquad
g_k(x)=
\frac{\widetilde g(\eta^kx)-a_k\eta^kx_n}{\eta^kA_k}.
\]
The rescaled operator
$X\mapsto\eta^kA_k^{-1}\widetilde F(\eta^{-k}A_kX)$
has the same ellipticity constants as $\widetilde F$
and vanishes at zero.
Thus the constants $\eta$ and $\delta$ in Lemma~\ref{k1}
remain the same at every scale.

We now check the smallness assumptions of Lemma~\ref{k1}.
By \eqref{s1} and the induction hypothesis,
\[
\|u_k\|_{L^\infty((\Omega_k)_1)}
=
\frac{
\|\widetilde u-a_kx_n\|_{L^\infty(\widetilde\Omega_{\eta^k})}
}{\eta^kA_k}
\leq1.
\]
By \eqref{o1} and \eqref{lip3},
\[
A_k\geq\omega(\eta^k)\geq\eta^{k/(1+\gamma)}.
\]
Hence
\[
\|f_k\|_{L^\infty((\Omega_k)_1)}
=\eta^kA_k^{-(1+\gamma)}
\|\widetilde f\|_{L^\infty(\widetilde\Omega_{\eta^k})}
\leq
\|\widetilde f\|_{L^\infty(\widetilde\Omega_1)}
\leq\delta.
\]
To estimate $g_k$, we first show that
\begin{equation}\label{ak}
|a_k|\leq2\bar C.
\end{equation}
Indeed, by \eqref{lip3}, for every $N\geq1$,
\[
\sum_{j=0}^N A_j
\leq1+\sum_{j=1}^N\omega(\eta^j)
+\eta^\alpha\sum_{j=0}^{N-1}A_j.
\]
The monotonicity of $\omega$ and \eqref{a3} give
\[
\sum_{j=1}^N\omega(\eta^j)
\leq
\omega(\eta)+\frac{1}{1-\eta}
\sum_{j=2}^N\int_{\eta^j}^{\eta^{j-1}}\frac{\omega(t)}{t}\,dt
\leq
\omega(\eta)+\frac{1}{1-\eta}
\int_0^\eta\frac{\omega(t)}{t}\,dt
\leq\frac12.
\]
Since $\eta^\alpha\leq1/4$, it follows that
\[
\sum_{j=0}^N A_j
\leq\frac{3}{2(1-\eta^\alpha)}
\leq2,
\qquad
\sum_{j=0}^\infty A_j\leq2.
\]
For $k\geq1$, \eqref{lip2} and the above estimate yield
\[
|a_k|
\leq\sum_{i=1}^k|a_i-a_{i-1}|
\leq\bar C\sum_{j=0}^{k-1}A_j
\leq2\bar C.
\]
The same bound holds for $k=0$.
By \eqref{a3}, \eqref{ak}, and $A_k\geq\omega(\eta^k)$,
\begin{equation}\label{g1}
\|g_k\|_{L^\infty((\partial\Omega_k)_1)}
\leq
\frac{1}{\eta^kA_k}
\left(
\frac{\delta}{2}\eta^k\omega(\eta^k)
+2\bar C\frac{\delta}{4\bar C}\eta^k\omega(\eta^k)
\right)
=\delta\frac{\omega(\eta^k)}{A_k}
\leq\delta.
\end{equation}
The pointwise domain inclusions and \eqref{a3} also give
\begin{equation}\label{g2}
\mathop{\mathrm{osc}}_{B_1}\partial\Omega_k
\leq
\frac{1}{\eta^k}
\mathop{\mathrm{osc}}_{B_{\eta^k}}\partial\widetilde\Omega
\leq
\frac{\delta}{2\bar C}\omega(\eta^k)
\leq\delta.
\end{equation}

We can now apply Lemma~\ref{k1}.
There exists $\bar a\in\mathbb R$ such that
\[
|\bar a|\leq\bar C,
\qquad
\|u_k-\bar a x_n\|_{L^\infty((\Omega_k)_\eta)}
\leq\eta^{1+\alpha}.
\]
Define
\[
a_{k+1}=a_k+\bar a A_k.
\]
Then
\[
|a_{k+1}-a_k|
=|\bar a|A_k
\leq\bar C A_k,
\]
and, by \eqref{lip3},
\[
\|\widetilde u-a_{k+1}x_n\|
_{L^\infty(\widetilde\Omega_{\eta^{k+1}})}
=\eta^kA_k
\|u_k-\bar a x_n\|_{L^\infty((\Omega_k)_\eta)}
\leq\eta^{k+1}\eta^\alpha A_k
\leq\eta^{k+1}A_{k+1}.
\]
Thus \eqref{lip1}--\eqref{lip2} hold for $k+1$,
and the induction is complete.

\par\smallskip\noindent
\textbf{Step 3. We pass to the limit and return to the original variables.}
By \eqref{lip1}--\eqref{lip2},
\[
\|\widetilde u-a_kx_n\|_{L^\infty(\widetilde\Omega_{\eta^k})}
\leq\eta^kA_k,
\qquad
|a_{k+1}-a_k|\leq\bar C A_k
\quad(k\geq0).
\]
For $m>k$, we have
\[
|a_m-a_k|
\leq\sum_{j=k}^{m-1}|a_{j+1}-a_j|
\leq\bar C\sum_{j=k}^{m-1}A_j.
\]
Since $\sum_{j=0}^\infty A_j<\infty$, there exists a constant $a$ such that
\[a=\lim_{k\to\infty}a_k,\qquad
|a-a_k|\leq\bar C\sum_{j=k}^\infty A_j,
\qquad
|a|\leq2\bar C.
\]
For any $x\in\widetilde\Omega_\eta$, choose $k\geq1$ such that
\[
\eta^{k+1}\leq|x|<\eta^k.
\]
By \eqref{lip1} and the above estimates,
\[
\begin{aligned}
|\widetilde u(x)-ax_n|
\leq
|\widetilde u(x)-a_kx_n|+|a_k-a|\,|x_n|
\leq
\eta^kA_k+\bar C|x|\sum_{j=k}^\infty A_j
\leq
\eta^k\left(A_k+\bar C\sum_{j=k}^\infty A_j\right).
\end{aligned}
\]
Define $\widetilde\sigma(0)=0$ and
\[
\widetilde\sigma(\eta^k)
=A_{k-1}+\bar C\sum_{j=k-1}^\infty A_j
\quad(k\geq1).
\]
{Interpolate linearly on each interval
$[\eta^{k+1},\eta^k]$, $k\geq1$.}
Since $\omega(\eta^k)\leq A_{k-1}$, we derive from
\eqref{lip3} that $A_k\leq A_{k-1}$.
Thus $\{A_k\}$ is nonincreasing and summable,
so $\widetilde\sigma$ is a modulus of continuity.
For $\eta^{k+1}\leq|x|<\eta^k$, we obtain
\[
|\widetilde u(x)-ax_n|
\leq\eta^k\widetilde\sigma(\eta^{k+1})
\leq\eta^{-1}|x|\widetilde\sigma(|x|).
\]
This proves \eqref{star}.

We now return to the original variables.
Recall that
\[
\widetilde u(x)
=\frac{u(rx)-g(0)-rDg(0)\cdot x}{M}.
\]
Define
\[
L(y)=g(0)+Dg(0)\cdot y+\frac{M}{r}a y_n.
\]
For every $y\in\Omega_{r\eta}$, we have
$y/r\in\widetilde\Omega_\eta$.
Hence
\[
\begin{aligned}
|u(y)-L(y)|
=M\left|
\widetilde u\!\left(\frac{y}{r}\right)
-a\frac{y_n}{r}
\right|
\leq
\frac{CM}{r}|y|\widetilde\sigma\!\left(\frac{|y|}{r}\right).
\end{aligned}
\]
Moreover,
\[
|DL|
=\left|Dg(0)+\frac{M}{r}a e_n\right|
\leq |Dg(0)|+\frac{2\bar C M}{r}.
\]
Set $r_0=r\eta$ and define
\[
\sigma(t)=\widetilde\sigma(t/r)
\quad\text{for }0\leq t\leq r_0,
\]
extending $\sigma$ constantly for $t\geq r_0$.
Since $|Dg(0)|\leq M$, these estimates give
\eqref{l12}--\eqref{l13}.
\qed

\begin{re}\label{dini}
The pointwise $C^{1,\mathrm{Dini}}$ assumptions on
$\partial\Omega$ and $g$ give the boundary bounds used in
\eqref{g1}--\eqref{g2}.
The Dini condition also ensures that
$\sum_{k=0}^{\infty}A_k<\infty$.
Thus the affine slopes converge, and the function $\sigma$
constructed above satisfies $\sigma(r)\to0$ as $r\to0$.
Without the Dini condition, this argument does not guarantee
these conclusions.
\end{re}

\bigskip

\noindent\textbf{Proof of Theorem~\ref{t1.2}.}

\par\smallskip\noindent
\textbf{Step 1. We normalize the problem.}
Let $\alpha_0=\min\{\alpha,1/(1+\gamma)\}$.
Choose $\eta$ as in Lemma~\ref{k1} with exponent $\alpha_0$,
taking $\eta$ small enough that $\eta^{\alpha_0}\leq1/2$.
Then fix the corresponding $0<\delta\leq1$.
Set
\[
M=\|u\|_{L^\infty(\Omega_1)}
  +\|f\|_{L^\infty(\Omega_1)}^{1/(1+\gamma)}
  +\|g\|_{C^{1,\alpha}(0)}.
\]
The case $M=0$ is immediate.
Assume $M>0$ and choose $0<r\leq1$ such that
\[
r^{\gamma+2}\leq\delta,\qquad
r^{1+\alpha}\leq\frac{\delta}{2},\qquad
r^\alpha\|\partial\Omega\|_{C^{1,\alpha}(0)}
\leq\frac{\delta}{4\bar C}.
\]
Define
\[
\widetilde\Omega=\Omega/r,\qquad
\widetilde u(x)
=\frac{u(rx)-g(0)-rDg(0)\cdot x}{M},
\qquad
p=\frac{rDg(0)}{M}.
\]
Then $\widetilde u$ satisfies
\[
\begin{cases}
|D\widetilde u+p|^\gamma
\widetilde F(D^2\widetilde u)=\widetilde f
&\text{in }\widetilde\Omega_1,\\
\widetilde u=\widetilde g
&\text{on }(\partial\widetilde\Omega)_1,
\end{cases}
\]
where
\[
\widetilde F(X)=\frac{r^2}{M}
F\!\left(\frac{M}{r^2}X\right),
\qquad
\widetilde f(x)
=\frac{r^{\gamma+2}}{M^{\gamma+1}}f(rx), \qquad \widetilde g(x)
=\frac{g(rx)-g(0)-rDg(0)\cdot x}{M}.
\]
The operator $\widetilde F$ has the same ellipticity constants
as $F$ and satisfies $\widetilde F(0)=0$.
By the choices of $M$ and $r$, the normalized quantities satisfy
\begin{equation}\label{a3'}
\begin{gathered}
\|\widetilde u\|_{L^\infty(\widetilde\Omega_1)}\leq1,
\quad
\|\widetilde f\|_{L^\infty(\widetilde\Omega_1)}\leq\delta\qquad \mathrm{and}\qquad
|\widetilde g(x)|\leq\frac{\delta}{2}|x|^{1+\alpha},
\quad
|x_n|\leq\frac{\delta}{4\bar C}|x|^{1+\alpha}
\quad\text{on }(\partial\widetilde\Omega)_1.
\end{gathered}
\end{equation}
We will prove that there exists $a\in\mathbb R$ such that
\begin{equation}\label{lxm}
|\widetilde u(x)-ax_n|
\leq C|x|^{1+\alpha_0}
\quad\text{for }x\in\widetilde\Omega_\eta,
\qquad |a|\leq C,
\end{equation}
where $C$ depends only on
$n,\lambda,\Lambda,\gamma$, and $\alpha_0$.

\par\smallskip\noindent
\textbf{Step 2. We construct affine approximations by induction.}
We construct a sequence of constants $\{a_k\}_{k\geq0}$ such that
\begin{equation}\label{3.10}
\begin{gathered}
\|\widetilde u-a_kx_n\|_{L^\infty(\widetilde\Omega_{\eta^k})}
\leq\eta^{k(1+\alpha_0)}
\qquad(k\geq0),\\
|a_k-a_{k-1}|
\leq\bar C\eta^{(k-1)\alpha_0}
\qquad(k\geq1),
\end{gathered}
\end{equation}
where $\bar C$ is as in Lemma~\ref{l3.1}.

For $k=0$, set $a_0=0$.
The first estimate in \eqref{3.10} follows from
$\|\widetilde u\|_{L^\infty(\widetilde\Omega_1)}\leq1$.
Assume that \eqref{3.10} holds for all indices up to $k$.
We prove it for $k+1$.
Set
\begin{equation}\label{s1-holder}
\Omega_k=\widetilde\Omega/\eta^k,
\qquad
u_k(x)=
\frac{\widetilde u(\eta^kx)-a_k\eta^kx_n}
{\eta^{k(1+\alpha_0)}},
\qquad x\in(\Omega_k)_1.
\end{equation}
Then $u_k$ satisfies
\[
\left\{
\begin{aligned}
&|Du_k+\eta^{-k\alpha_0}(p+a_ke_n)|^\gamma
\eta^{k(1-\alpha_0)}
\widetilde F\!\left(
\eta^{k(\alpha_0-1)}D^2u_k
\right)=f_k
&&\text{in }(\Omega_k)_1,\\
&u_k=g_k
&&\text{on }(\partial\Omega_k)_1,
\end{aligned}
\right.
\]
where
\[
f_k(x)=
\eta^{k[1-\alpha_0(1+\gamma)]}\widetilde f(\eta^kx),
\qquad
g_k(x)=
\frac{\widetilde g(\eta^kx)-a_k\eta^kx_n}
{\eta^{k(1+\alpha_0)}}.
\]
The rescaled operator
$
X\longmapsto
\eta^{k(1-\alpha_0)}
\widetilde F\!\left(\eta^{k(\alpha_0-1)}X\right)
$
has the same ellipticity constants as $\widetilde F$
and vanishes at zero.
Thus the constants $\eta$ and $\delta$ in Lemma~\ref{k1}
remain the same at every scale.

We now check the smallness assumptions of Lemma~\ref{k1}.
By the induction hypothesis,
\[
\|u_k\|_{L^\infty((\Omega_k)_1)}
=
\frac{
\|\widetilde u-a_kx_n\|_{L^\infty(\widetilde\Omega_{\eta^k})}
}{
\eta^{k(1+\alpha_0)}
}
\leq1.
\]
Since $\alpha_0\leq1/(1+\gamma)$,
\[
\begin{aligned}
\|f_k\|_{L^\infty((\Omega_k)_1)}
=\eta^{k[1-\alpha_0(1+\gamma)]}
\|\widetilde f\|_{L^\infty(\widetilde\Omega_{\eta^k})}
\leq
\|\widetilde f\|_{L^\infty(\widetilde\Omega_1)}
\leq\delta.
\end{aligned}
\]
To estimate $g_k$, we first bound $a_k$.
The induction hypothesis gives
\[
\begin{aligned}
|a_k|
\leq\sum_{i=1}^k|a_i-a_{i-1}|
\leq\bar C\sum_{j=0}^{k-1}\eta^{j\alpha_0}
\leq\frac{\bar C}{1-\eta^{\alpha_0}}
\leq2\bar C,
\end{aligned}
\]
where we used $\eta^{\alpha_0}\leq1/2$.
By \eqref{a3'}, the above bound, and $\alpha_0\leq\alpha$,
\begin{equation}\label{g1'}
\begin{aligned}
\|g_k\|_{L^\infty((\partial\Omega_k)_1)}
\leq
\frac{1}{\eta^{k(1+\alpha_0)}}
\left(
\frac{\delta}{2}\eta^{k(1+\alpha)}
+
2\bar C\frac{\delta}{4\bar C}\eta^{k(1+\alpha)}
\right)
=\delta\eta^{k(\alpha-\alpha_0)}
\leq\delta.
\end{aligned}
\end{equation}
The pointwise domain inclusions and \eqref{a3'} also give
\begin{equation}\label{g2'}
\mathop{\mathrm{osc}}_{B_1}\partial\Omega_k
\leq
\frac{1}{\eta^k}
\mathop{\mathrm{osc}}_{B_{\eta^k}}\partial\widetilde\Omega
\leq
\frac{\delta}{2\bar C}\eta^{k\alpha}
\leq\delta.
\end{equation}

We can now apply Lemma~\ref{k1}.
There exists $\bar a\in\mathbb R$ such that
\[
|\bar a|\leq\bar C,
\qquad
\|u_k-\bar a x_n\|_{L^\infty((\Omega_k)_\eta)}
\leq\eta^{1+\alpha_0}.
\]
Define
\[
a_{k+1}=a_k+\bar a\eta^{k\alpha_0}.
\]
Then
\[
|a_{k+1}-a_k|
=|\bar a|\eta^{k\alpha_0}
\leq\bar C\eta^{k\alpha_0},
\]
and
\[
\begin{aligned}
\|\widetilde u-a_{k+1}x_n\|
_{L^\infty(\widetilde\Omega_{\eta^{k+1}})}
=\eta^{k(1+\alpha_0)}
\|u_k-\bar a x_n\|_{L^\infty((\Omega_k)_\eta)}
\leq\eta^{(k+1)(1+\alpha_0)}.
\end{aligned}
\]
Thus \eqref{3.10} holds for $k+1$, and the induction is complete.

\par\smallskip\noindent
\textbf{Step 3. We pass to the limit and return to the original variables.}
By \eqref{3.10}, we have
\[
\|\widetilde u-a_kx_n\|_{L^\infty(\widetilde\Omega_{\eta^k})}
\leq\eta^{k(1+\alpha_0)},
\qquad
|a_{k+1}-a_k|\leq\bar C\eta^{k\alpha_0}
\quad(k\geq0).
\]
For $m>k$, it follows that
\[
\begin{aligned}
|a_m-a_k|
\leq\sum_{j=k}^{m-1}|a_{j+1}-a_j|
\leq\bar C\sum_{j=k}^{m-1}\eta^{j\alpha_0}
\leq\frac{\bar C}{1-\eta^{\alpha_0}}\eta^{k\alpha_0}
\leq2\bar C\eta^{k\alpha_0},
\end{aligned}
\]
where we used $\eta^{\alpha_0}\leq1/2$.
Thus $\{a_k\}$ is a Cauchy sequence and there exists a constant $a$ such that  
\[a=\lim_{k\to\infty}a_k,\qquad
|a-a_k|\leq2\bar C\eta^{k\alpha_0},
\qquad
|a|\leq2\bar C.
\]
For any $x\in\widetilde\Omega_\eta$, choose $k\geq1$ such that
\[
\eta^{k+1}\leq|x|<\eta^k.
\]
By \eqref{3.10} and the above estimates,
\[
\begin{aligned}
|\widetilde u(x)-ax_n|
\leq
|\widetilde u(x)-a_kx_n|+|a_k-a|\,|x_n|
\leq
\eta^{k(1+\alpha_0)}
+2\bar C\eta^{k\alpha_0}|x|
\leq C\eta^{k(1+\alpha_0)}
\leq C|x|^{1+\alpha_0}.
\end{aligned}
\]
This proves \eqref{lxm}.
Here, $C$ depends only on $n$, $\lambda$, $\Lambda$,
$\gamma$, and $\alpha_0$.

We now return to the original variables.
Recall that
\[
\widetilde u(x)
=\frac{u(rx)-g(0)-rDg(0)\cdot x}{M}.
\]
Define
\[
L(y)=g(0)+Dg(0)\cdot y+\frac{M}{r}a y_n.
\]
For every $y\in\Omega_{r\eta}$, we have
$y/r\in\widetilde\Omega_\eta$.
Hence
\[
\begin{aligned}
|u(y)-L(y)|
=M\left|
\widetilde u\!\left(\frac{y}{r}\right)
-a\frac{y_n}{r}
\right|
\leq
\frac{CM}{r^{1+\alpha_0}}|y|^{1+\alpha_0}.
\end{aligned}
\]
Moreover,
\[
|DL|
=\left|Dg(0)+\frac{M}{r}a e_n\right|
\leq |Dg(0)|+\frac{2\bar C M}{r}.
\]
These estimates give \eqref{l1}--\eqref{l2}
with the stated dependence of the constants.
\qed

\begin{re}\label{alpha}
As in the preceding remark, the pointwise $C^{1,\alpha}$
assumptions on $\partial\Omega$ and $g$ give the bounds
for the rescaled boundary data and boundary oscillation
in \eqref{g1'}--\eqref{g2'}.
These bounds allow the induction to continue at every scale.
\end{re}

\section{A Counterexample}\label{sec5}

\begingroup
\emergencystretch=3em
 In this section, we construct a counterexample showing that  $W^{2,\delta}$ regularity does not hold {for} $|Du|^\gamma F(D^2u)=f$ with $\gamma>0$. More precisely, {no uniform interior $W^{2,\delta}$ estimate follows from boundedness of $u$ and $f$ alone}.
\par\endgroup

\medskip
\noindent\textbf{Proof of Theorem \ref{t1.4}.}
\enlargethispage{3\baselineskip}
{Fix $\gamma>0$ and set $\epsilon=\gamma/(1+\gamma)\in(0,1)$.}
Let
\[
{u_M(x)}=M^{\epsilon-2}\sin Mx_1,\qquad x\in B_1,\quad M>1.
\]
Clearly, $\|{u_M}\|_{L^\infty(B_1)}\le1$. A direct calculation gives
\[
{\partial_1u_M}=M^{\epsilon-1}\cos Mx_1,
\qquad {\partial_{11}u_M}=-M^\epsilon\sin Mx_1.
\]
{Since $\gamma(\epsilon-1)+\epsilon=0$,}
\[
|D{u_M}|^{{\gamma}}\Delta{u_M}
={f_M}:=-|\cos Mx_1|^{{\gamma}}\sin Mx_1,
\qquad \|{f_M}\|_{L^\infty(B_1)}\le1.
\]
{
For $\delta>0$, integration over the cylinder
$(-1/4,1/4)\times B'_{1/4}\subset B_{1/2}$ yields
\[
\begin{aligned}
\int_{B_{1/2}}\|D^2u_M\|^\delta\,dx
\ge c_n M^{\epsilon\delta}\int_{-1/4}^{1/4}|\sin Mt|^\delta\,dt
\ge c_{n,\delta}M^{\epsilon\delta}\longrightarrow\infty.
\end{aligned}
\]
Consequently, for any $\delta>0$, $||D^2 u||_{L^\delta(B_1)}$ {cannot} be controlled by $||u||_{L^{\infty}(B_1)}$ and $||f||_{L^{\infty}(B_1)}$, {showing} that $W^{2,\delta}$ regularity {does not hold}.
}
\qed





\Needspace{8\baselineskip}
\begin{center}
	\textbf{Acknowledgements}	
\end{center}
\quad\quad The first author is supported by {the} is supported by National Science Foundation of China (No.12471202) and Shaanxi Fundamental Science Research Project for Mathematics and Physics (No.25JSY047). The second author is {supported by the} National Natural Science Foundation of China (No.12401257) and the Natural Science Foundation of Shandong Province, China (No.ZR2024QA045).
\begin{samepage}
\begin{center}
	\textbf{Data Availability Statement}	
\end{center}
\quad\quad Data will be made available on request.
\begin{center}
	\textbf{Declarations}	
\end{center}
\quad\quad The authors declare {that} they have no financial interests.
\end{samepage}
\bigskip

\bibliographystyle{elsarticle-num}


\enlargethispage{2\baselineskip}

\end{document}